\documentclass[11pt]{amsart}
\usepackage{amssymb}
\theoremstyle{plain}
\newtheorem{theorem}{Theorem}
\newtheorem{corollary}{Corollary}

\newtheorem{proposition}{Proposition}

\theoremstyle{definition}

\theoremstyle{remark}

\numberwithin{equation}{section}
\def\K#1#2{\displaystyle{\mathop K\limits_{#1}^{#2}}}

\setbox0=\hbox{$+$}
\newdimen\plusheight
\plusheight=\ht0
\def\+{\;\lower\plusheight\hbox{$+$}\;}

\setbox0=\hbox{$-$}
\newdimen\minusheight
\minusheight=\ht0
\def\-{\;\lower\minusheight\hbox{$-$}\;}

\setbox0=\hbox{$\cdots$}
\newdimen\cdotsheight
\cdotsheight=\plusheight
\def\cds{\lower\cdotsheight\hbox{$\cdots$}}

\begin{document}
\title[ Extensions and Contractions of Continued Fractions ]
 {Ramanujan and   Extensions and Contractions of Continued Fractions }
\author{J. Mc Laughlin}
\address{Mathematics Department\\
 Trinity College\\
300 Summit Street, Hartford, CT 06106-3100}
\email{james.mclaughlin@trincoll.edu}
\author{ Nancy J. Wyshinski}
\address{Mathematics Department\\
       Trinity College\\
        300 Summit Street, Hartford, CT 06106-3100}
\email{nancy.wyshinski@trincoll.edu}

\date{January 29, 2003}
\begin{abstract}
If a continued fraction $K_{n=1}^{\infty} a_{n}/b_{n}$ is known to converge but its limit is
not easy to determine,
 it may be easier to use an  extension   of  $K _{n=1}^{\infty}a_{n}/b_{n}$ to
find the limit. By an extension of $K_{n=1}^{\infty} a_{n}/b_{n}$ we mean a continued fraction
$K_{n=1}^{\infty} c_{n}/d_{n}$ whose odd or even part is $K_{n=1}^{\infty} a_{n}/b_{n}$.
One can then possibly find the limit in one of three ways:

(i) Prove the extension converges and find its limit;

(ii) Prove the extension converges and find the limit of the other contraction
(for example, the odd part, if $K _{n=1}^{\infty}a_{n}/b_{n}$ is the even part);

(ii) Find the limit of the other contraction and show that the odd and
even parts of the extension tend to the same limit.

We apply these ideas to derive new proofs of certain continued
fraction identities of Ramanujan and to prove a generalization of
an identity involving the Rogers-Ramanujan continued fraction,
  which was conjectured by Blecksmith and Brillhart.
\end{abstract}
\maketitle

%           Introduction

\section{Introduction}
The methods used by the great Indian mathematician, Srinivasa Ramanujan, to obtain many
of his fascinating results remain a mystery. In this paper we
describe some simple ideas concerning extensions and
 contractions of continued
fractions which may have
led Ramanujan to some of the elegant entries concerning  continued fractions
that  he made in his famous notebooks (See \cite{B89}, Chapter 12).

Suppose we are given a continued fraction
$d_{0}+K_{n=1}^{\infty}c_{n}/d_{n}$ whose limit is sought.
If the limit is difficult
to compute, it may be easier to work with one of  several
 extensions of the continued fraction which can easily
be written down.
Suppose, for example, that the even part of an extension gives the
original continued fraction.
One can  try to find the limit in one of three other ways:\\
(i) Prove the extension converges and find its limit;\\
(ii) Prove the extension converges and find the limit of the odd part;\\
(iii) Find the limit of the odd part and show that the odd and
even parts of the extension tend to the same limit (by showing that
the absolute value of the difference between consecutive approximants
of the extension tends to
$0$).

For (i) above,  an equivalence transformation of
the extended continued fraction may result in its approximants naturally taking on
a particulary simple form, so that its limit (and thus of its even part --
the original continued fraction) can be found
very easily (See, for example, the proofs of Entry 12 and Entry 13).

For (ii) and (iii), it may turn out that the limit
of the odd part of the extension can be computed almost trivially.
If one can then show that either the extension converges \emph{or}
that the even and odd parts of the extension tend to the same limit, then one knows that the limit
of the original continued fraction is the same as the limit of the odd part of the extension.
The advantage here is that it is usually much easier to show convergence than to determine the
actual limit of a continued fraction.

This is our reason for believing that Ramanujan may have used extensions
and contractions in the discovery of some of his results -- in several entries
the odd part of an extension of the continued fraction being considered by Ramanujan
can be shown to converge to give Ramanujan's claimed limit  almost trivially.

We illustrate the principles involved by giving new proofs of some of Ramanujan's
continued fraction
identities  found in Chapter 12 of the second notebook.

We also use these methods to give a generalization of an identity
involving the famous Rogers-Ramanujan continued fraction, first
conjectured by  Blecksmith and  Brillhart \cite{BB02}, and proved
by Berndt and Yee in \cite{BY02}.

\section{Extensions and Contractions of Continued Fractions}

We start with the concepts of extensions and
contractions of continued fractions. Before coming to details, we borrow some notation from \cite{LW92}
(page 83).\footnote{The authors mention
in \cite{LW92} that this idea also
goes back to Seidel \cite{S55} and that Lagrange had some special cases already in 1774 \cite{L74}
 and 1776 \cite{L76}.}
A continued fraction $d_{0}+K_{n=1}^{\infty}c_{n}/d_{n}$ is said to be a \emph{contraction} of the
continued fraction $b_{0}+K_{n=1}^{\infty}a_{n}/b_{n}$ if its classical approximants
$\{g_{n}\}$ form a subsequence of the classical approximants $\{f_{n}\}$
 of $b_{0}+K_{n=1}^{\infty}a_{n}/b_{n}$. In this case $b_{0}+K_{n=1}^{\infty}a_{n}/b_{n}$
is called an \emph{extension} of $d_{0}+K_{n=1}^{\infty}c_{n}/d_{n}$.

We call $d_{0}+K_{n=1}^{\infty}c_{n}/d_{n}$  a \emph{canonical contraction} of
 $b_{0}+K_{n=1}^{\infty}a_{n}/b_{n}$ if
\begin{align*}
&C_{k}=A_{n_{k}},& &D_{k}=B_{n_{k}}& &\text{ for } k=0,1,2,3,\ldots \, ,\phantom{asdasd}&
\end{align*}
where $C_{n}$, $D_{n}$, $A_{n}$ and $B_{n}$ are canonical numerators and denominators
of $d_{0}+K_{n=1}^{\infty}c_{n}/d_{n}$ and $b_{0}+K_{n=1}^{\infty}a_{n}/b_{n}$ respectively.

Here we use the standard notations
\[
\K{n = 1}{N}\displaystyle{\frac{a_{n}}{b_{n}}} \,:=
\cfrac{a_{1}}{b_{1} + \cfrac{a_{2}}{b_{2} + \cfrac{a_{3}}{b_{3} +
\ldots +  \cfrac{a_{N}^{}}{b_{N} }}}}
  \]
\[
\phantom{asd} =  \frac{a_{1}}{b_{1}+}\, \frac{a_{2}}{b_{2}+}\,
\frac{a_{3}}{b_{3}+}\, \dots  \frac{a_{N}}{b_{N}}.
\]
We write $\displaystyle{A_{N}/B_{N}}$ for the above finite continued
fraction written as a rational function of the variables
$a_{1},...,a_{N},b_{1},...,b_{N}$.
By $K_{n = 1}^{\infty}\displaystyle{a_{n}/b_{n}}$
 we  mean the limit
of the sequence  \{$\displaystyle{A_{n}/B_{n}}$\} as \,$n$\,goes to infinity,
if the limit exists. The ratio $A_{N}/B_{N}$ is called the $N$-th \emph{approximant}
of the continued fraction. It is elementary that the $A_{N}$ (the $N$th (canonical) \emph{numerator})
and $B_{N}$ (the $N$th (canonical) \emph{denominator})
satisfy the following
recurrence relations:
\begin{align}\label{E:recur}
A_{N}&=b_{N}A_{N-1}+a_{N}A_{N-2},\\
B_{N}&=b_{N}B_{N-1}+a_{N}B_{N-2} \notag.
\end{align}
It can also be easily shown that
\begin{align}\label{detform}
A_{N}B_{N-1}-A_{N-1}B_{N}&=(-1)^{n-1}\prod_{i=1}^{N}a_{i},\\
A_{N+1}B_{N-1}-A_{N-1}B_{N+1}&=(-1)^{n-1}b_{N+1}\prod_{i=1}^{N}a_{i}. \notag
\end{align}

From \cite{LW92} (page 83) we have the following theorem:
\begin{theorem}\label{T:t1}
The canonical contraction of $b_{0}+K_{n=1}^{\infty}a_{n}/b_{n}$ with
\begin{align*}
&C_{k}=A_{2k}& &D_{k}=B_{2k}& &\text{ for } k=0,1,2,3,\ldots \, ,&
\end{align*}
exists if and only if $b_{2k} \not = 0 for K=1,2,3,\ldots$, and in this case is given by
\begin{equation}\label{E:evcf}
b_{0}
+
\frac{b_{2}a_{1}}{b_{2}b_{1}+a_{2}}
\-
\frac{a_{2}a_{3}b_{4}/b_{2}}{a_{4}+b_{3}b_{4}+a_{3}b_{4}/b_{2}}
\-
\frac{a_{4}a_{5}b_{6}/b_{4}}{a_{6}+b_{5}b_{6}+a_{5}b_{6}/b_{4}}
\+
\cds .
\end{equation}
\end{theorem}
The continued fraction \eqref{E:evcf} is called the \emph{even} part of $b_{0}+K_{n=1}^{\infty}a_{n}/b_{n}$.

We give some simple corollaries to this theorem, which we will use
later.
\begin{corollary}\label{corcf1}
The even part of
\begin{equation}\label{cf1}
d_{0}
+
\frac{c_{1}}{d_{1}-c_{2}}
\+
\frac{c_{2}}{1}
\+
\frac{-1}{d_{2}-c_{3}+1}
\+
\frac{c_{3}}{1}
\+
\frac{-1}{d_{3}-c_{4}+1}
\+
\frac{c_{4}}{1}
\+
\cds
\end{equation}
is
\begin{equation}\label{cf2}
d_{0}
\+
\frac{c_{1}}{d_{1}}
\+
\frac{c_{2}}{d_{2}}
\+
\frac{c_{3}}{d_{3}}
\+
\frac{c_{4}}{d_{4}}
\+
\cds .
\end{equation}
\end{corollary}
\begin{proof}
In Theorem \ref{T:t1}, set $b_{0}$ = $d_{0}$, $a_{1}=c_{1}$, $b_{1}=d_{1}-c_{2}$ and, for $k\geq 1$,
$a_{2k}=c_{k+1}$, $a_{2k+1}=-1$, $b_{2k}=1$ and $b_{2k+1}=d_{k+1}-c_{k+2}+1$.
\end{proof}

\begin{corollary}\label{corcf2}
The even
part of
\begin{equation}\label{cf3}
d_{0}
+
\frac{c_{1}}{d_{1}-1}
\+
\frac{-1}{1}
\+
\frac{c_{2}}{d_{2}-c_{2}+1}
\+
\frac{-1}{1}
\+
\frac{c_{3}}{d_{3}-c_{3}+1}
\+
\frac{-1}{1}
\+
\cds
\end{equation}
is
\begin{equation}\label{cf2a}
d_{0}
\+
\frac{c_{1}}{d_{1}}
\+
\frac{c_{2}}{d_{2}}
\+
\frac{c_{3}}{d_{3}}
\+
\frac{c_{4}}{d_{4}}
\+
\cds .
\end{equation}
\end{corollary}
\begin{proof}
In Theorem \ref{T:t1}, set  $b_{0}$ = $d_{0}$, $a_{1}=c_{1}$, $b_{1}=d_{1}+1$ and, for $k\geq 1$,
$a_{2k}=-1$, $a_{2k+1}=c_{k+1}$, $b_{2k}=1$ and $b_{2k+1}=d_{k+1}-c_{k+1}+1$.
\end{proof}

\begin{corollary}\label{C:cor3a}
The even part
of
\begin{equation}\label{cf4}
b_{0}
+
\frac{a_{1}}{b_{1}}
\+
\frac{a_{2}}{1}
\+
\frac{a_{3}}{0}
\+
\frac{a_{4}}{1}
\+
\frac{a_{5}}{0}
\+
\frac{a_{6}}{1}
\+
\frac{a_{7}}{0}
\+
\cds
\end{equation}
is
\begin{equation}\label{cf5}
b_{0}
+
\frac{a_{1}}{b_{1}+a_{2}}
\-
\frac{a_{2}a_{3}}{a_{4}+a_{3}}
\-
\frac{a_{4}a_{5}}{a_{6}+a_{5}}
\-
\frac{a_{6}a_{7}}{a_{8}+a_{7}}
\-
\cds .
\end{equation}
\end{corollary}
\begin{proof}
In Theorem \ref{T:t1}, set $b_{2k}=1$ and $b_{2k+1}=0$,  for $k \geq 1$.
\end{proof}

From \cite{LW92} (page 85) we also have:
\begin{theorem}\label{odcf}
The canonical contraction of $b_{0}+K_{n=1}^{\infty}a_{n}/b_{n}$ with
$C_{0}=A_{1}/B_{1}$
\begin{align*}
&C_{k}=A_{2k+1}& &D_{k}=B_{2k+1}& &\text{ for } k=1,2,3,\ldots \, ,&
\end{align*}
exists if and only if $b_{2k+1} \not = 0 for K=0,1,2,3,\ldots$, and in this case is given by
\begin{multline}\label{E:odcf}
\frac{b_{0}b_{1}+a_{1}}{b_{1}}
-
\frac{a_{1}a_{2}b_{3}/b_{1}}{b_{1}(a_{3}+b_{2}b_{3})+a_{2}b_{3}}
\-
\frac{a_{3}a_{4}b_{5}b_{1}/b_{3}}{a_{5}+b_{4}b_{5}+a_{4}b_{5}/b_{3}}\\
\-
\frac{a_{5}a_{6}b_{7}/b_{5}}{a_{7}+b_{6}b_{7}+a_{6}b_{7}/b_{5}}
\-
\frac{a_{7}a_{8}b_{9}/b_{7}}{a_{9}+b_{8}b_{9}+a_{8}b_{9}/b_{7}}
\+
\cds .
\end{multline}
\end{theorem}
The continued fraction \eqref{E:odcf} is called the \emph{odd} part of $b_{0}+K_{n=1}^{\infty}a_{n}/b_{n}$.

We will also make use of the following corollary to Theorem \ref{odcf}.
\begin{corollary}\label{corcf7}
The odd part of the continued fraction
\begin{equation}\label{cf7}
\frac{c_{1}}{1}
\-
\frac{c_{2}}{1}
\+
\frac{c_{2}}{1}
\-
\frac{c_{3}}{1}
\+
\frac{c_{3}}{1}
\-
\frac{c_{4}}{1}
\+
\frac{c_{4}}{1}
\-
\cds
\end{equation}
is
\begin{equation}\label{cf8}
c_{1}
+
\frac{c_{1}c_{2}}{1}
\+
\frac{c_{2}c_{3}}{1}
\+
\frac{c_{3}c_{4}}{1}
\+
\cds .
\end{equation}
\end{corollary}
\begin{proof}
In Theorem \ref{odcf}, set $b_{0}=0$, $a_{1}=c_{1}$, and, for $k\geq 1$, $b_{i}=1$, $a_{2i}=-c_{i}$
and $a_{2i+1}=c_{i}$.
\end{proof}
We will not explicitly
compute the odd parts of the continued fractions at \eqref{cf1},  \eqref{cf3} and \eqref{cf4} at this point.

We  also  give a new extension/contraction proof of Daniel Bernoulli's
transformation of a sequence into a continued fraction
\cite{B75} (see, for example, \cite{K63}, pp. 11--12).
\begin{proposition}
Let $\{K_{0},K_{1}, K_{2},\ldots\}$ be a sequence of complex numbers such that $K_{i}\not = K_{i-1}$,
for $i=1,2,\ldots$.

Then  $\{K_{0},K_{1}, K_{2},\ldots\}$  is the sequence of approximants of the
continued fraction
\begin{multline}\label{ber1}
K_{0}+\frac{K_{1}-K_{0}}{1}
\+
\frac{K_{1}-K_{2}}{K_{2}-K_{0}}
\+
\frac{(K_{1}-K_{0})(K_{2}-K_{3})}
            {K_{3}-K_{1}}
\+\\
\ldots
\+
\frac{(K_{n-2}-K_{n-3})(K_{n-1}-K_{n})}
            {K_{n}-K_{n-2}}
\+
\ldots
.
\end{multline}
\end{proposition}
\begin{proof}
We use the fact that
\[
\frac{1}{a}
\+
\frac{1}{0}
\+
\frac{1}{b}
\+
\frac{1}{c}
=
\frac{1}{a+b}
\+
\frac{1}{c}.
\]
Then
\begin{align*}
K_{0} &+
\frac{1}{0}
\+
\frac{1}{K_{1}-K_{0}}
\+
\frac{1}{0}
\+
\frac{1}{K_{2}-K_{1}}
\+
\cds
\+
\frac{1}{0}
\+
\frac{1}{K_{n}-K_{n-1}}\\
&=K_{0} +
\frac{1}{0}
\+
\frac{1}{\sum_{i=1}^{n}K_{i}-K_{i-1}}=K_{n}.
\end{align*}
On the other hand, by Theorem \ref{T:t1},  the even part of the above continued fraction
is
\begin{multline*}
K_{0}
+
\frac{K_{1}-K_{0}}{1}
\-
\frac{(K_{2}-K_{1})/(K_{1}-K_{0})}{1+(K_{2}-K_{1})/(K_{1}-K_{0})}
\-
\frac{(K_{3}-K_{2})/(K_{2}-K_{1})}{1+(K_{3}-K_{2})/(K_{2}-K_{1}) }\\
\-
\cds
\-
\frac{(K_{n}-K_{n-1})/(K_{n-1}-K_{n-2})}{1+ (K_{n}-K_{n-1})/(K_{n-1}-K_{n-2})}
\end{multline*}
%&\phantom{as}\\
\begin{multline*}
=K_{0}+\frac{K_{1}-K_{0}}{1}
\+
\frac{K_{1}-K_{2}}{K_{2}-K_{0}}
\+
\frac{(K_{1}-K_{0})(K_{2}-K_{3})}
            {K_{3}-K_{1}}
\+\\
\ldots
\+
\frac{(K_{n-2}-K_{n-3})(K_{n-1}-K_{n})}
            {K_{n}-K_{n-2}}.
\end{multline*}
\end{proof}
If we let $K_{n}=\sum_{i=0}^{n}a_{i}$,
we of course get
Euler's transformation of a series into a continued fraction:
\begin{equation}\label{eul1}
\sum_{i=0}^{n}a_{i} = a_{0}+\frac{a_{1}}{1}
\+
\frac{-a_{2}}{a_{2}+a_{1}}
\+
\frac{-a_{1}a_{3}}{a_{3}+a_{2}}
\+
\cds
\frac{-a_{n-2}a_{n}}{a_{n}+a_{n-1}}.
\end{equation}

\section{Some Continued Fractions from Chapter 12 of Ramanujan's Second Notebook}

In each of the following example, $b_{0}+K_{n=1}^{\infty}a_{n}/b_{n}$ will mean
the extended continued fraction under consideration and
$\{A_{n}\}$ and $\{B_{n}\}$ will denote its sequences numerators and
denominators respectively.

We will frequently make use of the following
important theorem of Worpitzky (see \cite{LW92}, pp. 35--36).

\begin{theorem}(Worpitzky)
 Let the continued fraction $K_{n=1}^{\infty}a_{n}/1$ be such that
$|a_{n}|\leq 1/4$ for $n \geq 1$. Then$K_{n=1}^{\infty}a_{n}/1$ converges.
 All approximants of the continued fraction lie in the disc $|w|<1/2$ and the value of the
continued fraction is in the disk $|w|\leq1/2$.\end{theorem}

We now illustrate the methods involved by giving new proofs
of several continued fraction identities due to Ramanujan.

\textbf{Entry 7.} ( \cite{B89}, page 112)  \emph{If $x$ is not a negative integer, then }
\begin{equation}\label{en7}
1=
\frac{x+1}{x}
\+
\frac{x+2}{x+1}
\+
\frac{x+3}{x+2}
\+
\cds .
\end{equation}
 We prove a generalization of Entry 7.

\textbf{Entry 7a.}  \emph{ Let $\{y_{i}\}_{i=1}^{\infty}$ be any sequence of complex numbers such that}
\begin{align*}
&(i)\,y_{i} \not = -1, \,\,i=1,2,3,\ldots ,\\
&(ii)\,\lim_{n \to \infty} \prod_{i=1}^{n}|1+y_{i}|= \infty, \\
&(iii)\left| \frac{y_{i}+1}{y_{i-1}\,y_{i}} \right| \leq \frac{1}{4},  \text{ \emph{for all }}
i\geq N_{0}, \text{ some } N_{0}.
\end{align*}
\emph{Then}
\begin{equation}\label{en7a}
1=
\frac{y_{1}+1}{y_{1}}
\+
\frac{y_{2}+1}{y_{2}}
\+
\frac{y_{3}+1}{y_{3}}
\+
\cds
\end{equation}
\begin{proof}
It is sufficient to assume $N_{0}=2$ (if not, one
 can prove the result for the  tail that
begins with $y_{N_{0}}$ in the numerator and the
 continued fraction will then collapse from
the bottom up to give the result).

After a similarity transformation,
the left side of  Equation \ref{en7a} becomes
\[
\frac{(y_{1}+1)/y_{1}}{1}
\+
\frac{(y_{2}+1)/(y_{1}y_{2})}{1}
\+
\frac{(y_{3}+1)/(y_{2}y_{3})}{1}
\+
\cds
\]
 and Worpitzky's theorem   gives that this continued fraction, and thus the left side of
\eqref{en7a}  converges.

From Corollary \ref{corcf2} it can be see that the right side of Equation \ref{en7a} is
the even part of
\begin{equation}\label{cf7aex}
\frac{y_{1}+1}{y_{1}+1}
\+
\frac{-1}{1}
\+
\frac{y_{2}+1}{0}
\+
\frac{-1}{1}
\+
\frac{y_{3}+1}{0}
\+
\frac{-1}{1}
\+
\cds .
\end{equation}
 Thus the even part of \eqref{cf7aex} converges and,
from Equation \ref{detform},
\begin{equation*}
0=\lim_{i \to \infty}\left|\frac{A_{2i+2}}{B_{2i+2}}-\frac{A_{2i}}{B_{2i}}\right|=
\lim_{i \to \infty}\left|\frac{A_{2i+2}B_{2i}-A_{2i}B_{2i+2}}{B_{2i+2}B_{2i}}\right|=
\lim_{i \to \infty}\frac{\prod_{j=1}^{i+1}|y_{j}+1|}{|B_{2i}B_{2i+2}|}
\end{equation*}
Condition (ii) above then gives that
\begin{equation}\label{Blim}
\lim_{i \to \infty}|B_{2i}B_{2i+2}|=\infty.
\end{equation}
From the recurrence relations at \eqref{E:recur} and the fact that $b_{2i+1}=0$, for   $i=1,2,\ldots$,
one has that
\[
A_{2i+1}=(y_{i+1}+1)A_{2i-1}= \dots = \prod_{j=1}^{i+1}(y_{j}+1).
\]
 Similarly,
$B_{2i+1}= \prod_{j=1}^{i+1}(y_{j}+1)$ and so each odd-numbered
approximant is identically 1.
From Equation \ref{detform} above it also follows that
\[
|A_{2i+1}B_{2i}-A_{2i}B_{2i+1}|=\prod_{j=1}^{i+1}|y_{j}+1|=|B_{2i+1}|.
\]
Thus
\begin{equation}
\left| 1-\frac{A_{2i}}{B_{2i}}\right|=
\left| \frac{A_{2i+1}}{B_{2i+1}}-\frac{A_{2i}}{B_{2i}}\right|=
\left| \frac{A_{2i+1}B_{2i}-A_{2i}B_{2i+1}}{B_{2i}B_{2i+1}}\right|=
\frac{1}{|B_{2i}|}.
\end{equation}
Since the limit of the left side exists, it follows
that $\lim_{i \to \infty}|B_{2i}|$ exists
and Equation \ref{Blim} gives that this limit is $\infty$.
Thus $\lim_{i \to \infty}A_{2i}/B_{2i}=1$ and the result follows.
\end{proof}

\vspace{5pt}

For the next example, we will show that the extended
continued fraction converges, so that the even and odd parts have the same limit.
We will give two different proofs to better illustrate the methods.
 One will use the following theorem of Lange \cite{L66} (see \cite{JT80}, page 124):

\begin{theorem}\label{T:lange}
The continued fraction $K(c_{n}^{2}/1)$ converges to a finite value provided that
\begin{equation}\label{langeq}
|c_{2n-1}\pm i \alpha| \leq \rho, \,\,\,\,|c_{2n}\pm i (1+\alpha)| \geq \rho, \,\,\,\,n=1,2,3,\ldots ,
\end{equation}
where $\alpha$ is a complex number and $\alpha$
and $\rho$ satisfy the inequality
\begin{equation}\label{langeq2}
|\alpha|<\rho<|\alpha+1|.
\end{equation}
The convergence is uniform with respect to the regions defined by \eqref{langeq}.
\end{theorem}
(We have changed the notation in the above theorem slightly to avoid
conflict with existing notation.)

The other proof will use the following theorem,
due to Wall \cite{W56} (see \cite{JT80}, page 127):

\begin{theorem}\label{T:wall}
Let $\{f_{n}\}$ be the sequence of approximants of a continued fraction $K(a_{n}/1)$. Assume
that there exist positive numbers $M$, $L$, $n_{0}$ and a
subsequence $\{m_{k}\}$ of the positive integers such that
\begin{align}\label{walleq1}
&\text{\phantom{ad}} |f_{n}|<M, & & \text{\phantom{asdaghjksd}}n=n_{0},n_{0}+1,n_{0}+2,\ldots ,&
\end{align}
and
\begin{align}\label{walleq2}
&|a_{m_{k}}|<L, \text{\phantom{asdasd}}& & k=1,2,3,\ldots .\text{\phantom{asdasd}}&
\end{align}
Further assume that the odd (even) part of $K(a_{n}/1)$ converges to a finite value $v$. Then there
exists a subsequence of the even (odd) part which converges to $v$.
\end{theorem}

Remark:  An obvious implication of this theorem is that if,
in addition, the odd and even parts both converge,
then they converge to the same limit and the continued fraction converges to this limit.

\textbf{Entry 9.} (\cite{B89}, page 114) \emph{Let $a$ and $x$ be complex numbers such that either $x\not = - k \, a$
for $k \in \{1,2,\ldots\}$ and $a \not = 0$, or that $a = 0$ and $|x|>1$. Then }
\begin{equation}\label{en9}
\frac{x+a+1}{x+1}
=
\frac{x+a}{x-1}
\+
\frac{x+2a}{x+a-1}
\+
\frac{x+3a}{x+2a-1}
\+
\cds .
\end{equation}

We will not consider the case $a=0$, since the right side is periodic for $a=0$
and the result follows from a general theorem for periodic continued fractions.

Neither will we consider the case $x=-1$, since this case
("$\infty=\infty$") follows from the case $x \not = -1$ by
considering the tail beginning with $-1+3a$ in the numerator.

FIRST PROOF.  We will prove Entry 9 for $a \not \in (-\infty,\,0)$.
From Corollary \ref{corcf2}, the right side of Equation \ref{en9} is the even part of
\begin{equation}\label{enex}
\frac{x+a}{x} \+ \frac{-1}{1} \+ \frac{x+2a}{-a} \+ \frac{-1}{1}
\+ \frac{x+3a}{-a} \+ \frac{-1}{1} \+ \frac{x+4a}{-a} \+ \cds
\end{equation}
From Theorem \ref{odcf}, the odd part of this latter continued fraction is
{\allowdisplaybreaks
\begin{align*}
&\frac{x+a}{x} +
\frac{(x+a)(-a)/x}{x(x+a)+a}
\+
\frac{(x+2a)x}{(x+2a)-1}\\
&\+
\frac{(x+3a)}{(x+3a)-1}
\+
\frac{(x+4a)}{(x+4a)-1}
\+
\cds\\
&\phantom{as}\\
&=\frac{x+a}{x} +
\frac{(x+a)(-a)/x}{x(x+a)+a+x}
&\phantom{as}\\
&=\frac{x+a+1}{x+1}
\end{align*}
}
The second last equality follows from Entry 7a applied to the tail of the continued fraction.
Thus the odd part of the extension converges to the left side of Equation \ref{en9}.
Our first proof   that the extension itself converges uses Theorem \ref{T:lange}.
The continued fraction at \eqref{enex} is equivalent to the following continued fraction:
\begin{multline}\label{enexeq}
\frac{(x+a)/x}{1}
\+
\frac{-1/x}{1}
\+
\frac{-x/a-2}{1}
\+
\frac{1/a}{1}
\+\\
\frac{-x/a-3}{1}
\+
\frac{1/a}{1}
\+
\frac{-x/a-4}{1}
\+
\cds
\end{multline}
We consider a tail of this continued fraction
\begin{equation}\label{enext}
\frac{1/a}{1}
\+
\frac{-x/a-m}{1}
\+
\frac{1/a}{1}
\+
\frac{-x/a-m-1}{1}
\+
\frac{1/a}{1}
\+
\frac{-x/a-m-2}{1}
\+
\cds ,
\end{equation}
where $m$ will depend on $a$ and $x$ and will be determined later.
If the tail converges, then the continued fraction converges and its limit is $(x+a+1)/(x+1)$, since the
odd part converges to this limit. Denote the continued fraction at \eqref{enext} by
$K_{k=1}^{\infty}c_{k}^{2}/1$, so that
\begin{align}
&c_{2k-1}^{2}=\frac{1}{a},&  &c_{2k}^{2}=-\frac{x}{a}-m-k+1.&
\end{align}
Note that if $\alpha$ and $\rho$ can be found such that inequality
\eqref{langeq2} and the first inequality in \eqref{langeq} can be
satisfied, then the second inequality in \eqref{langeq} will be
satisfied automatically for all $k$, provided $m$ is chosen large
enough. Let $\sqrt{1/a}=c+i\, d$, where $c>0$ (since $a \not \in
(-\infty,\,0)$).  Set
\begin{align}
&\alpha=\frac{c^2+d^2}{2}\left(1+i \frac{d}{c}\right),&
&\rho=\sqrt{\frac{c^2+d^2}{4c^2}((c^2+d^2)^2+4c^2)}.&
\end{align}
Then $|\alpha|<\rho<|\alpha+1|=\sqrt{\rho^2+1}$ and
$|c_{2k-1}+ i \alpha| =|c_{2k-1}- i \alpha|=\rho$. Provided $m$ is chosen large enough
so that $|-x/a-m-k+1 \pm i(1+\alpha)| \geq \rho$, for $k=1,2,\ldots ,$ the conditions of
Theorem \ref{T:lange} are satisfied, the tail converges to a finite value and the extended
continued fraction converges and Entry 9 follows for $a \not \in (-\infty, 0)$.
\begin{flushright}
$\Box$
\end{flushright}
SECOND PROOF. We will apply Theorem \ref{T:wall} to the continued fraction at
\eqref{enexeq}.
 Without loss of generality,
we can assume that
\[
\left|\frac{x+ja}{(x+(j-2)a-1)(x+(j-1)a-1)}\right|\leq \frac{1}{4}
\]
 holds for
$j\geq 2$, since this holds for all $j$ sufficiently large and if Entry 9
 holds for a tail of the continued fraction, ie.,
\begin{multline}\label{en9a}
\frac{x+na+1}{x+(n-1)a+1} = \frac{x+n a}{x+(n-1)a-1}
\+\\
\frac{x+(n+1)a}{x+na-1}
\+
\frac{x+(n+2)a}{x+(n+1)a-1}
\+
\cds ,
\end{multline}
for some integer $n$, then Entry 9 is proved,
since the continued fraction will then collapse from the bottom up to give the result.
 Thus Worpitzky's Theorem gives that the even part of \ref{enex}, and
thus the even part of \ref{enexeq}, converges to a finite value
and we already know that the odd part converges to the left side of Equation \ref{en9}.
Thus, since the odd and even parts both converge to finite values,
 there exist an $n_{0}$ and an $M$ such that Equation \ref{walleq1} is satisfied.
Further, from \eqref{enexeq}, it is clear that in Equation \ref{walleq2},
we can take $\{m_{k}\}$ to be the even integers
and $L$ to be $\max \{|1/a|+1, |1/x|+1\}$.
(The case $x=0$ is not a problem since we can consider the tail of the
continued fraction beginning with the third partial numerator.)
The conditions of Theorem \ref{T:wall} are satisfied and, by the remark
following it, the continued fraction at \eqref{enex} converges and Entry 9 follows.
\begin{flushright}
$\Box$
\end{flushright}

\textbf{Entry 10} (\cite{B89}, page 116)  \emph{If $n$ is a positive integer, then}
\begin{equation}\label{en10}
n=
\frac{1}{1-n}
\+
\frac{2}{2-n}
\+
\frac{3}{3-n}
\+
\cds
\+
\frac{n}{0}
\+
\frac{n+1}{1}
\+
\frac{n+2}{2}
\+
\cds .
\end{equation}

\begin{proof}
The proof is by induction on $n$. If $n=1$, the left side of \eqref{en10} is
\begin{equation*}
n=
\frac{1}{1-1}
\+
\frac{2}{1}
\+
\frac{3}{2}
\+
\frac{4}{3}
\+
\cds
=
\frac{1}{1-1+1} =1,
\end{equation*}
by Entry 7.  Suppose Entry 10 is true for $n=1,2, \ldots , m-1$. From Corollary \ref{corcf2},
\begin{equation}\label{em10}
\frac{1}{1-m}
\+
\frac{2}{2-m}
\+
\frac{3}{3-m}
\+
\cds
\+
\frac{m}{0}
\+
\frac{m+1}{1}
\+
\frac{m+2}{2}
\+
\cds
\end{equation}
is the even part of
\begin{equation}\label{enex10}
\frac{1}{2-m}
\+
\frac{-1}{1}
\+
\frac{2}{1-m}
\+
\frac{-1}{1}
\+
\frac{3}{1-m}
\+
\cds .
\end{equation}
From Theorem \ref{odcf}, the  odd part of this latter continued fraction is
\begin{align}\label{en10od}
&\frac{1}{2-m}
+
\frac{(1-m)/(2-m)}{(2-m)(3-m)-(1-m)}
\+
\frac{2(2-m)}{(4-m)-1}\\
&\+
\frac{3}{(5-m)-1}
\+
\frac{4}{(6-m)-1}
\+
\cds \notag \\
&\phantom{as}\notag \\
&=\frac{1}{2-m}
+
\frac{(1-m)/(2-m)}{(2-m)(3-m)-(1-m)}
\+
\frac{(2-m)2}{2-(m-1)}\notag\\
&\+
\frac{3}{3-(m-1)}
\+
\frac{4}{4-(m-1)}
\+
\cds
\notag\\
&\phantom{as}\notag \\
&=\frac{1}{2-m}
+
\frac{(1-m)/(2-m)}
{  (2-m)(3-m)-(1-m) +(2-m)
\left(
\frac{1}{m-1}-(2-m)
\right)
} \notag\\
&=m.\notag
\end{align}
The next to last step comes from applying the induction step to
the continued fraction in \eqref{en10}, when $n=m-1$.

Next, we will use Theorem \ref{T:lange} to show that a tail of the continued fraction
at \eqref{enex10} converges and thus that the continued fraction itself converges to $m$,
since the odd part
equals $m$. This continued fraction is equivalent
to
\begin{multline}\label{en10exeq}
\frac{1/(2-m)}{1}
\+
\frac{1/(m-2)}{1}
\+
\frac{2/(1-m)}{1}\\
\+
\frac{1/(m-1)}{1}
\+
\frac{3/(1-m)}{1}
\+
\frac{1/(m-1)}{1}
\+
\cds .
\end{multline}
We consider a tail of this last continued fraction:
\begin{multline}\label{en10exeqt}
\frac{1/(m-1)}{1}
\+
\frac{N/(1-m)}{1}
\+
\frac{1/(m-1)}{1}
\+
\frac{(N+1)/(1-m)}{1}
\+
\cds ,
\end{multline}
where $N$ depends on $m$ and will be chosen later. With the notation of Theorem \ref{T:lange},
let this continued fraction be denoted $K_{k=1}^{\infty}(c_{k}^{2}/1)$. Since $m\geq 2$,
$c_{2k-1}=\sqrt{1/(m-1)}$ is real and we chose the positive square root so that
$c_{2k-1}=:c>0$.  Let $c_{2k} = + i  \sqrt{(N+k-1)/(m-1)}$. We chose
 $\alpha = c^2/2$ and $\rho = \sqrt{c^4/4+c^2}$. Then
\[
|\alpha| <  \rho <  |\alpha+1|=\sqrt{\rho^2+1}.
\]
Further,
\[
|c_{2k-1}+i \, \alpha|=|c_{2k-1}-i \, \alpha| = \rho.
\]
For $N$ sufficiently large, $|c_{2k}\pm i (1+\alpha)| \geq \rho$.
The conditions of Theorem \ref{T:lange} are satisfied,  a tail of
\eqref{en10exeq} converges and, by the remark following the statement of this
theorem, \eqref{en10exeq} itself converges and Entry 10 follows.
\end{proof}

\textbf{Entry 12.}  (\cite{B89}, page 118) \emph{If $a \not  = 0$ and $x \not = - k\,a$,
where $k$ is a positive integer,}
\begin{equation}\label{en12}
1=
\frac{x+a}{a}
\+
\frac{(x+a)^2-a^2}{a}
\+
\frac{(x+2 a)^2-a^2}{a}
\+
\frac{(x+3 a)^2-a^2}{a}
\+
\cds .
\end{equation}

\begin{proof}
The left side of Equation \ref{en12} is
\begin{multline}\label{en12s}
\frac{x+a}{a}
\+
\frac{x(x+2 a)}{a}
\+
\frac{(x+a)(x+3 a)}{a}
\+\\
\cds
\+
\frac{(x+(k-1) a)(x+(k+1)a)}{a}
\+
\cds ,
\end{multline}
which, by Corollary \ref{C:cor3a},  is the even part of
\begin{multline}\label{enex12}
\frac{x+a}{x+a}
\+
\frac{-x}{1}
\+
\frac{x+2a}{0}
\+
\frac{-(x+a)}{1}
\+
\frac{x+3a}{0}
\+
\frac{-(x+2a)}{1}
\+
\cds .
\end{multline}
This continued fraction is equivalent to the following continued fraction:
\begin{multline}\label{enex12den}
\frac{1}{1}
\+
\frac{1}{-1-a/x}
\+
\frac{1}{0}
\+
\frac{1}{1+2 a/x}
\+
\frac{1}{0}
\+
\frac{1}{-1-3 a/x}\\
\+
\frac{1}{0}
\+
\frac{1}{1+4 a/x}
\+
\frac{1}{0}
\+
\frac{1}{-1-5 a/x}
\+
\frac{1}{0}
\+
\frac{1}{1+6a/x}
\+
\cds
\end{multline}
By the similar reasoning to that used in the proof of Entry 7a,  each odd-indexed
 approximant of
\eqref{enex12den} is identically $1$. To calculate the even-indexed approximants,
we make repeated use the identity
\begin{equation}\label{zden}
\frac{1}{a}
\+
\frac{1}{0}
\+
\frac{1}{b}
\+
\frac{1}{c}
=
\frac{1}{a+b}
\+
\frac{1}{c}
\end{equation}
to simplify the approximant.
One easily checks that
\begin{align}
%\label{enex12denevp}
&\frac{1}{1}
\+
\frac{1}{-1-a/x}
\+
\frac{1}{0}
\+
\frac{1}{1+2 a/x}
\+
\frac{1}{0}
\+
\cds
\+
\frac{1}{0}
\+
\frac{1}{1+2k a/x} \notag \\
&=\frac{1}{1}
\+
\frac{1}{(-1-a/x)+(1+2a/x)+(-1-3a/x)+\dots+(1+2k a/x)} \notag
\\
&=\frac{1}{1}
\+
\frac{1}{ka/x}. \notag
\end{align}
Similarly, it is easy to check that
\begin{align}\label{enex12denodm}
&\frac{1}{1}
\+
\frac{1}{-1-a/x}
\+
\frac{1}{0}
\+
\frac{1}{1+2 a/x}
\+
\frac{1}{0}
\+
\cds
\+
\frac{1}{0}
\+
\frac{1}{-1-(2k-1) a/x} \notag \\
&=\frac{1}{1}
\+
\frac{1}{(-1-a/x)+(1+2a/x)+(-1-3a/x)+\dots+(-1-(2k-1) a/x)} \notag
\\
&=\frac{1}{1}
\+
\frac{1}{-1-ka/x}. \notag
\end{align}
Upon letting $k \to \infty$, one has that the even-indexed tend to 1 also and Entry 12
is proved.
\end{proof}

Before coming to Entry 13, we state the following theorem of Hill \cite{H08} (see:
\cite{AAR99}, page 63)
\begin{theorem}\label{T:hill}
Let $s_{n}$ denote the $n$th partial sum of $_{2}F_{1}(a,b;c;1)$. For
\text{Re}$(c-a-b)<0$,
\begin{equation}\label{hyp1}
s_{n}\sim \frac{\Gamma (c)n^{a+b-c}}{\Gamma (a) \Gamma (b)},
\end{equation}
and for $c=a+b$,
\begin{equation}\label{hyp2}
s_{n} \sim \frac{\Gamma (c)\log n}{\Gamma (a) \Gamma (b)}.
\end{equation}
\end{theorem}
Here
\[
_{2}F_{1}(a,b;c;x)=\sum_{n=0}^{\infty}\frac{(a)_{n}(b)_{n}}{(c)_{n}n!}x^{n},
\]
where  $(d)_{n}=d(d+1)\ldots(d+n-1)$ for $n>0$ and  $(d)_{0}=1$.

\textbf{Entry 13.} (\cite{B89}, page 119)
\emph{Let $a$, $b$ and $d$ be complex numbers such that either $d \not = 0$,
$b \not = - k d$, where $k$ is a non-negative integer, and} Re(($a-b)/d)  <0$,
\footnote{ Entry 13, as written in \cite{B89}, page 119, reads "Re($(a-b)/d)  >0$". This
inequality should be reversed. The following example is an indication of this fact
(a complete proof that the reversed inequality is the
correct one is found in the proof of Entry 13 above): Let $a=2$ and $b=d=1$, so that
Re$((a-b)/d)>0$. However the left side of
\eqref{en13} is
\begin{align*}
&\frac{2.1}{4}
\+
\frac{-3.2}{6}
\+
\frac{-4.3}{8}
\+
\frac{-5.4}{10}
\+
\frac{-6.5}{12}
\+
\cds \\
=&\frac{1}{2}
\+
\frac{-1}{2}
\+
\frac{-1}{2}
\+
\frac{-1}{2}
\+
\frac{-1}{2}
\+
\cds .
\end{align*}
This continued fraction has the sequence of approximants $\{n/(n+1)\}$ and
therefore converges to $1$ (=$b$) and not $2$ (=$a$).
 }
\emph{or $d\not = 0$ and $a=b$, or $d=0$ and $|a|<|b|$.
Then}
\begin{equation}\label{en13}
a
=
\frac{a b}{a+b+d}
\-
\frac{(a+d)(b+d)}{a+b+3 d}
\-
\frac{(a+2d)(b+2d)}{a+b+5d}
\-
\cds
\end{equation}
\begin{proof}
We make the further assumption
$a+kd \not =0$, for $k$ a non-negative integer.
By Corollary \ref{C:cor3a},
the right side of Equation \ref{en13} is the even part of the continued fraction
\begin{multline}\label{en13ex}
\frac{ab}{b}
\+
\frac{a+d}{1}
\+
\frac{b+d}{0}
\+
\frac{a+2d}{1}
\+
\frac{b+2d}{0}\\
\+
\cds
\frac{a+kd}{1}
\+
\frac{b+kd}{0}
\+
\cds .
\end{multline}
This latter continued fraction is equivalent to
\begin{multline}\label{en13exden}
\frac{1}{1/a}
\+
\frac{1}{\frac{ab}{a+d}}
\+
\frac{1}{0}
\+
\frac{1}{\frac{ab(b+d)}{(a+d)(a+2 d)}}
\+
\frac{1}{0}
\+ \\
\frac{1}{\frac{ab(b+d)(b+2d)}{(a+d)(a+2 d)(a+3d)}}
\+
\frac{1}{0}
\+
\frac{1}{\frac{ab(b+d)(b+2d)(b+3d)}{(a+d)(a+2 d)(a+3d)(a+4d)}}
\+
\frac{1}{0}
\+
\cds .
\end{multline}
By similar reasoning to that used in Example7a, each odd-indexed approximant
is identically equal to $a$. We now consider the even-indexed approximants, treating
each of the three cases in the statement of Entry 13 in turn. Note
 that, since $d \not =0$ (as in the first case),
\begin{align}\label{shiffac}
&\phantom{=}\frac{ab(b+d)\ldots(b+(k-1)d)}{(a+d)(a+2 d)\ldots(a+kd)}
=a\frac{(b/d)_{k}}{(a/d+1)_{k}}.
\end{align}

Let $\{f_{k}\}$ denote the sequence of approximants for the continued fraction
at \eqref{en13exden}.
By using the same collapsing technique as was used in
the proof of Example 12,
\begin{align}\label{evapprox}
f_{2k}=\frac{1}{1/a}
\+
\frac{1}{ a\displaystyle{\sum_{i=1}^{k}\frac{(b/d)_{i}}{(a/d+1)_{i}}}}
=\frac{1}{1/a}
\+
\frac{1}{ -a+a\displaystyle{\sum_{i=0}^{k}\frac{(b/d)_{i}}{(a/d+1)_{i}}}}
\end{align}
We note that $\displaystyle{\sum_{i=0}^{k}\frac{(b/d)_{i}}{(a/d+1)_{i}}}$
is the $k$th partial sum of $_{2}F_{1}(1,b/d;a/d+1;1)$.
Since Re$(a/d+1-b/d-1)$ $=$ Re$((a-b)/d)<0$, we have by Theorem \ref{T:hill} that
\begin{equation}
\displaystyle{\sum_{i=0}^{k}\frac{(b/d)_{i}}{(a/d+1)_{i}}}
\sim
 \frac{\Gamma (a/d+1)k^{1+b/d-a/d-1}}{\Gamma (1) \Gamma (b/d)}
=\frac{\Gamma (a/d+1)k^{(b-a)/d}}{\Gamma (1) \Gamma (b/d)}.
\end{equation}
Thus, since Re$((b-a)/d)>0$,
\[
\lim_{k \to \infty}\sum_{i=0}^{k}\frac{(b/d)_{i}}{(a/d+1)_{i}}=\infty
\]
and from Equation \ref{evapprox},  $\lim_{k \to \infty }f_{2k}=a$.

Next, suppose $d \not = 0$ and $a=b$. After canceling common
factors in each denominator and collapsing the continued fraction
as before, we have  that
 \begin{align}\label{evapproxaeqb}
f_{2k}=\frac{1}{1/a}
\+
\frac{1}{ \displaystyle{\sum_{i=1}^{k}\frac{a^{2}}{a+id}}},
\end{align}
and once again it is clear that   $\lim_{k \to \infty }f_{2k}=a$.

Finally, if $d=0$ and $|a|<|b|$
\begin{align}\label{evapproxaleb}
f_{2k}=\frac{1}{1/a}
\+
\frac{1}{a \displaystyle{\sum_{i=1}^{k}\left(\frac{b}{a}\right)^{i}}},
\end{align}
and once again  $\lim_{k \to \infty }f_{2k}=a$.
Entry 13 is proved.
\end{proof}
Remark:  Interestingly, this proof, deriving from extending the
right side of \eqref{en13}, coincides at the finish with
Jacobsen's proof  \cite{J89},  which uses a theorem, due to her
\cite{J86} and Waadeland \cite{Wa84}, on tails of continued
fractions.   Both proofs eventually  rely on Hill's result from
Theorem \ref{T:hill} applied to the same $_{2}F_{1}$ function.

\section{An Extension of the Rogers-Ramanujan Continued Fraction}
For $|q|<1$, let
\begin{equation}\label{rrcf}
R(q):=
1
+
\frac{q}{1}
\+
\frac{q^2}{1}
\+
\frac{q^3}{1}
\+
\cds ,
\end{equation}
the famous Rogers-Ramanujan continued fraction. In \cite{BB02}, J.
Brillhart and R. Blecksmith conjectured that
\begin{equation}
R(q)=
\frac{1}{1}
\-
\frac{q}{1}
\+
\frac{q}{1}
\-
\frac{q}{1}
\+
\frac{q}{1}
\-
\frac{q^2}{1}
\+
\frac{q^2}{1}
\-
\frac{q^2}{1}
\+
\frac{q^2}{1}
\-
\cds .
\end{equation}
This conjecture was proved by Berndt and Yee in \cite{BY02}. We
generalize this result as follows. For $q$, $\alpha \in
\mathbb{C}$, let $q^{\alpha}$ be defined as usual by
\[
q^{\alpha}=e^{\alpha \log q},
\]
where $\log q$ is the principal logarithm of $q$.
\begin{proposition}
Let $q$, $\alpha \in \mathbb{C}$, with $|q|<1$. Then
\begin{multline}\label{BBgen}
R(q)=1-q^{\alpha}
+
\frac{q^{\alpha}}{1}
\-
\frac{q^{1-\alpha}}{1}
\+
\frac{q^{1-\alpha}}{1}
\-
\frac{q^{1+\alpha}}{1}
\+
\frac{q^{1+\alpha}}{1}\\
\-
\frac{q^{2-\alpha}}{1}
\+
\frac{q^{2-\alpha}}{1}
\-
\frac{q^{2+\alpha}}{1}
\+
\frac{q^{2+\alpha}}{1}
\-
\cds .
\end{multline}
\end{proposition}
Remark: The conjecture of Blecksmith and Brillhart is the $\alpha=0$ case of this proposition.
\begin{proof}
In Corollary \ref{corcf7}, let $c_{1}=q^{\alpha}$ and, for $k \geq 1$, let
$c_{2k}=q^{k-\alpha}$ and $c_{2k+1}=q^{k+\alpha}$. This gives that the odd part of
\begin{equation*}
\frac{q^{\alpha}}{1}
\-
\frac{q^{1-\alpha}}{1}
\+
\frac{q^{1-\alpha}}{1}
\-
\frac{q^{1+\alpha}}{1}
\+
\frac{q^{1+\alpha}}{1}
\-
\frac{q^{2-\alpha}}{1}
\+
\frac{q^{2-\alpha}}{1}
\-
\frac{q^{2+\alpha}}{1}
\+
\frac{q^{2+\alpha}}{1}
\-
\cds
\end{equation*}
is
\[
q^{\alpha}+
\frac{q}{1}
\+
\frac{q^2}{1}
\+
\frac{q^3}{1}
\+
\cds .
\]
Since a tail of the left side of \eqref{BBgen} converges (by Worpitzky's Theorem) and its odd part
converges to $R(q)$, this proves the result.
\end{proof}
We also have the following corollary.
\begin{corollary}
Let $q$, $\alpha \in \mathbb{C}$, with $|q|<1$. Then
\begin{multline}\label{BBgenev}
R(q)=1-q^{\alpha}
+
\frac{q^{\alpha}}{1-q^{1-\alpha}}
\+
\frac{q^{2(1-\alpha)}}{1+q^{1-\alpha}-q^{1+\alpha}}
\+
\frac{q^{2(1+\alpha)}}{1+q^{1+\alpha}-q^{2-\alpha}}\\
\+
\frac{q^{2(2-\alpha)}}{1+q^{2-\alpha}-q^{2+\alpha}}
\+
\frac{q^{2(2+\alpha)}}{1+q^{2+\alpha}-q^{3-\alpha}}
\+
\cds
\end{multline}
\end{corollary}
\begin{proof}
By Theorem \ref{T:t1}, the left side of Equation \ref{BBgenev}
 is the even part of the continued fraction at \eqref{BBgen}.
\end{proof}

\allowdisplaybreaks{}

\end{document}